\documentclass[draft,12pt]{elsarticle}

\usepackage{amssymb}
\usepackage{amsthm}
\usepackage{amsmath,bbm,dsfont, color}
\usepackage{xfrac}

\journal{arXiv}

\begin{document}
\renewcommand{\theequation}{\thesection.\arabic{equation}}
\overfullrule5pt
\newcommand*{\id}{\mathds{1}}
\newcommand*{\zahlen}{\mathbb{Z}}
\newcommand*{\en}{\mathbb{N}}
\newcommand*{\er}{\mathbb{R}}
\newcommand{\esssup}{\mathop{\text{ess\:sup}}}
\newcommand{\essinf}{\mathop{\text{ess\:inf}}}
\newcommand{\osc}{\mathop{\text{osc}}}
\newcommand{\Inf}{\mathop{\text{inf}}}
\newcommand{\htimes}{\mathop{\text{\large$�$}}}
\newcommand{\hexists}{\mathop{\text{\LARGE$\exists$}}}
\newcommand{\hforall}{\mathop{\text{\LARGE$\forall$}}}
\newtheorem{theo}{Theorem}[section]                    
\newtheorem{lem}[theo]{Lemma}
\newtheorem{cor}[theo]{Corollary}
\newtheorem{prop}{Proposition}
\newtheorem{mydef}{Definition}[section]
\newtheorem{rem}{Remark}[section]      
\newcommand{\Sym}{{\rm Sym}}

\def\divv{{\rm div}}
\def\B{\mathbb{B}}
\def\D{\mathbb{D}}
\def\F{\mathbb{F}}
\def\N{\mathbb{N}}
\def\O{\mathbb{O}}
\def\R{\mathbb{R}}
\def\Kor{{\rm Kor}}
\def\pod#1{\mathop{#1}\limits}
\def\diagin{-\hskip-11.0truept\intop}
\def\diagint{{\raise-.1pt\hbox{--}\hskip-7.9pt\intop}}
\def\diagintop{\mathop{\mathchoice
{{\diagin}}%
{{\diagint}}%
{{\diagint}}%
{{\diagint}}%
}\limits}

\begin{frontmatter}

 \title{Boundary De Giorgi-Ladyzhenskaya classes and their application to regularity of swirl of Navier-Stokes}

\author{Jan Burczak 
} \ead{jb@impan.pl} \address{Institute of Mathematics, Polish Academy of Sciences, \'Sniadeckich 8, 00-950 Warsaw.}






\begin{abstract}
The embeddings theorem of space-boundary-type DeGiorgi-Ladyzhenskaya parabolic classes into H\"older spaces is presented, which is useful for regularity considerations for parabolic boundary value problems. Additionaly, the application of this theory to Navier-Stokes's swirl is presented. 
\end{abstract}

\begin{keyword}
DeGiorgi classes, swirl of Navier-Stokes, regularity of parabolic systems

\end{keyword}




\end{frontmatter}

\section{Introduction}
We present an unified treatment of embeddings of boundary-type DeGiorgi-Ladyzhenskaya parabolic  classes into H\"older spaces. This result serves the  regularity studies certain PDEs. Therefore we restrict ourselves to the case of space boundary and do not consider time-boundary, as in the class of PDEs which can be tackled by this theory, the local-in-time smoothness is standard. Generally we follow ideas of \cite{lsu}, where the case of boundary regularity is briefly mentioned. Here we provide clear and complete proofs and improve the original result qualitatively by obtaining better H\"older exponents, which is done in spirit of \cite{zaj}. Finally, the application of this theory to Navier-Stokes's swirl is presented. 
\section{Notation and preliminary results } We work with a following geometric objects 
\begin{itemize}
\item $\Omega_T$ denoting space-time cylinder $\Omega \times [-T; 0]$ with a domain $\Omega$ as its base,
\item $\Gamma \subset \partial \Omega \times (-T, 0)$ is the open part of the space boundary of  $\Omega_T$, in which vicinity we are interested in boundary regularity (in the case of Dirichlet data we need to have certain regularity of boundary data on $\Gamma$),
\item $Q(\rho, \tau)$ be, for a fixed point $(x_0, t_0) \in \Gamma$, a boundary cylinder $ B_\rho (x_0) \times (t_0 - \tau; t_0)$, which is small enough to satisfy $\partial \Omega_T \cap Q(\rho, \tau) = \Gamma \cap Q(\rho, \tau) $.
\end{itemize}
We will use also a following notation:
\begin{equation}
|f|_{V (\Omega_T)} \equiv \sup_{t \in [-T, 0]} |f (t)|_{2, \Omega} + |\nabla f|_{2, \Omega_T} 
\end{equation}

\begin{equation}
V(\Omega_T)  \equiv  \{ f \in L^2 (\Omega_T): |f|_{V (\Omega_T)} < \infty \}
\end{equation}
where $|f|_{2, U} \equiv \int_U |f|^2$ and $\nabla$ means space gradient. Observe that here we assume that $|f(t)|_{2, \Omega}  < \infty$ for every $t$. Let
\begin{equation}
\osc_U f = max_U f - min_U f
\end{equation}
\begin{equation}
A^f_{k, \rho} (t)  \equiv  \{ x \in B_\rho (x_0) \cap \Omega : f (x,t) > k \}
\label{defA}
\end{equation}
\begin{equation}
\mu_{k, \rho, \tau}  \equiv  \int_{t_0 - \tau}^{t_0} \mu^{\frac{r}{q}} (A^f_{k, \rho} (t) ) dt
\end{equation}
\begin{equation}
f^{(k)}  \equiv  (f-k)^+
\end{equation}
We introduce now classes  $B_N, B_D$ dependent on further specified parameters. The former is useful for showing boundary regularity for Neumann problems, the latter for Dirichlet problems. Let us define formal inequality:
\begin{equation}
\begin{aligned}
&\int_{B_{\rho} (x_0) \cap \Omega} | w^{(k)} (x, t_0) \xi (x, t_0)|^2 dx + \int_{Q_{\rho, \tau} \cap \Omega_T}| \nabla w^{(k)} (x, t) \xi (x, t)|^2 dx dt \cr
&\le \int_{B_\rho (x_0) \cap \Omega} | w^{(k)} (x, t_0 - \tau) \xi (x, t_0 - \tau)|^2 dx + \cr
& \gamma \left[ \int_{Q_{\rho, \tau} \cap \Omega_T}  ( | \nabla \xi|^2 + \xi |\xi,_t|) | w^{(k)}|^2+\left(  \int_{t_0 - \tau}^{t_0} \left( \int_{A_{k, \rho} (t) } \xi (x, t) dx \right)^\frac{r}{q} dt \right)^{\frac{2 (1+ \kappa)}{r}}   \right] 
\end{aligned}
\label{defB}
\end{equation}
\begin{mydef}
$u \in B_N (\Omega_T, M, \gamma, r, \delta, \kappa)$ iff
\begin{enumerate}
\item[(i)] $u$ is a pointwisely  defined representative of a function in $V(\Omega_T) \cap L^\infty (\Omega_T)$ and $|u|_{\infty, \Omega_T} \le M$.
\item[(ii)] Inequality \eqref{defB} with $w \equiv \pm u$ and $\frac{1}{r} + \frac{n}{2q} = \frac{n}{4} $ holds for any $k \ge \esssup_{Q_{\rho, \tau}} w - \delta $ and $\xi \in C (Q_{\rho, \tau}), \; 0 \le \xi \le 1, \xi \equiv 0$ on $\partial B_\rho (x_0) \times (t_0 - \tau, t_0)$.

\end{enumerate}
\label{defBn}
\end{mydef}

\begin{mydef}
$u \in B_D (\Omega_T,  M,  \gamma, r, \delta, \kappa; \Gamma, c_\Gamma, \beta)$ iff
\begin{enumerate}
\item[(i)] $u$ is a pointwisely  defined representative of a function in $ V(\Omega_T) \cap L^\infty (\Omega_T)$ and $|u|_{\infty, \Omega_T} \le M$
\item[(ii)] Inequality \eqref{defB} with $w \equiv \pm u$ and $\frac{1}{r} + \frac{n}{2q} = \frac{n}{4} $ holds for any $k \ge \max ( \esssup_{Q_{\rho, \tau}}   w - \delta, \; \max_{\Gamma \cap Q_{\rho, \tau}} w) $ and $\xi \in C (Q_{\rho, \tau}), \; 0 \le \xi \le 1, \xi \equiv 0$ on $\partial B_\rho (x_0) \times (t_0 - \tau, t_0)$.
\item[(iii)]  for an open set $\Gamma \subset \partial \Omega \times (-T, 0)$ holds $ \osc_{Q (\rho, \rho^2) } u_{| \Gamma}  \le c_\Gamma \rho^\beta$
\end{enumerate}
\label{defBd}
\end{mydef}
\begin{rem}
It is important that $\xi $ does not have to vanish on $\Gamma$. 
\end{rem}
\begin{rem}
One can unessentialiy generalize definitions \ref{defBn}, \ref{defBd} demanding that \eqref{defB} holds merely for functions $ \xi$, which cutoff certain cylinders $Q_{\rho, \tau}$. 
\end{rem}
A quotation of a few well-known results ends this section.
\begin{lem}
For nonnegative $h \in W^{1,1} (B_\rho)$, vanishing on $U$ of positive Lebesgue measure, holds a following generalized Poincare inequality:
\begin{equation}
\int_U h \eta \le K_P \rho^n \frac{ \mu^\frac{1}{n} (U)}{\mu (U_0)} \int_{B_\rho} | \nabla h| \eta
\end{equation}
where $\eta \equiv \eta(|x|) \in [0,1]$ and $\eta_{|U_0} \equiv 1$, $K_P = 2^n\left(\frac{1}{n}+ \omega_n \right)$ .
\label{gp}
\end{lem}
\begin{lem}
Assume that $\Omega$ is convex. For nonnegative $h \in W^{1,1} (B_\rho \cap \Omega)$, vanishing on $U$ of positive Lebesgue measure, holds a following generalized Poincare inequality:
\begin{equation}
\int_U h \eta \le \tilde{K}_P  \frac{\rho^{n+1}}{\mu (U_0)} \int_{B_\rho \cap \Omega} | \nabla h| \eta
\end{equation}
where $\eta \equiv \eta(|x|) \in [0,1]$ and $\eta_{|U_0} \equiv 1$, $\tilde{K}_P = $ .
\label{gpn}
\end{lem}
{\color{red} Suggestion of} Proof can be found in \cite{lsu}, p.92.
\section{Results}
\noindent
The following conditions excluding cusps of $\Omega$ are needed for validity of results
\begin{equation}
\hexists_{\theta_0> 0, \rho_0 > 0} \hforall_{\rho \le \rho_0, (x, t) \in \Gamma} \mu ( B_\rho (x)  \cap \Omega^c) \ge \theta_0 \mu (B_\rho (x) ),
\label{ac}
\end{equation}
\begin{equation}
\hexists_{\theta_0> 0, \rho_0 > 0} \hforall_{\rho \le \rho_0, (x, t) \in \Gamma} \mu ( B_\rho (x)  \cap \Omega) \ge \theta_0 \mu (B_\rho (x) );
\label{acn}
\end{equation}
the former allows for a much simplification of the result concerning Dirichlet boundary case and is referred to as {\it the anti-outer-cusp condition} in what follows. The latter plays a role at the Neumann boundary case and is referred to as {\it the anti-inner-cusp condition}. 
\begin{theo}
Assume that  {\it the anti-outer-cusp condition}  \eqref{ac} holds. \\
Take $u \in B_D (\Omega_T,  M,  \gamma, r, \delta, \kappa; \Gamma, c_\Gamma, \beta)$ with 
\[
r > 2 \; \text{ for } \; n = 2 \quad \text{ and } \quad r \ge 2 \;  \text{ for }\; n > 2
\]
\[
M < \infty, \; \gamma >0,  \; \delta > 0, \; \kappa > 0,  \; c_\Gamma < \infty, \; \beta > 0
\]

then $u$ is H\"older continuous in vicinity of $\Gamma$. \\
More precisely: take any $\sigma \in (1,2], \; \theta \in (0, 1]$ and a boundary cylinder $Q( \tilde{\rho}_0, \theta \tilde{\rho}_0^2) $  where $\tilde{\rho}_0 \le \rho_0$, with $ \rho_0$ from \eqref{ac},  such that $\partial \Omega_T \cap Q( \tilde{\rho}_0, \theta \tilde{\rho}_0^2) = \Gamma \cap Q( \tilde{\rho}_0, \theta \tilde{\rho}_0^2)$. We have for $\rho \le \sigma^{-2} \tilde{\rho}_0$ 
\begin{equation}
 \osc_{Q(\rho, \theta {\rho}^2)  \cap (\Omega_T \cup \Gamma)} u \le C \rho^\alpha
\end{equation}
with 
\begin{equation}
 \alpha = \min \left(- \log_{\sigma^2} ( 1 - 2^{-s}), \beta, \frac{n \kappa}{2} \right) \qquad C = \max \left( (\sfrac{\sigma^2}{\tilde{\rho}_0})^\alpha \max \left( \osc_{Q_{\tilde{\rho}_0}} u, 2^s \sigma^{\frac{n \kappa}{2}}  \tilde{\rho}_0^{\frac{n \kappa}{2}} \right), c_\Gamma \right)
\end{equation}
and $s$ satisfying
\begin{equation}
s \ge 1+ \max \left[ \left \lceil log_2 \frac{2 M }{\delta}\right \rceil\  +  \theta \left(4 \left(\frac{1}{n}+ \omega_n \right)^2  \omega_n^{\frac{2}{n}} \gamma \frac{2^{3n+2+2 \max \left(1, \frac{1+ \kappa}{r} \right) }}{ \eta^2 (\sigma - 1) \theta_0^2} \right),  \; log_2 (2 c_\Gamma \sigma^\beta) \right]
\end{equation}
\label{thm}
\end{theo}

\begin{theo}
Assume that  {\it the anti-outer-cusp condition}  \eqref{acn} holds. \\
Take $u \in B_N (\Omega_T,  M,  \gamma, r, \delta, \kappa)$ with 
\[
r > 2 \; \text{ for } \; n = 2 \quad \text{ and } \quad r \ge 2 \;  \text{ for }\; n > 2
\]
\[
M < \infty, \; \gamma >0,  \; \delta > 0, \; \kappa > 0
\]

then $u$ is H\"older continuous in vicinity of $\Gamma$. \\
More precisely: take any $\sigma \in (1,2]$,
\begin{equation}
\theta \le \min \left(1, \frac{\theta_0}{2304 \gamma}, \; \left(  \frac{\omega^{-\frac{2(1+ \kappa) - q}{q}}_n}{128 \gamma} \right)^{\frac{r}{2 (1+ \kappa)}} \right) 
\end{equation}

 and a boundary cylinder $Q( \tilde{\rho}_0, \theta \tilde{\rho}_0^2) $  where $\tilde{\rho}_0 \le \rho_0$, with $\rho_0$ from \eqref{acn}, such that $\partial \Omega_T \cap Q( \tilde{\rho}_0, \theta \tilde{\rho}_0^2) = \Gamma \cap Q( \tilde{\rho}_0, \theta \tilde{\rho}_0^2)$. We have for $\rho \le \sigma^{-2} \tilde{\rho}_0$ 
\begin{equation}
 \osc_{Q(\rho, \theta {\rho}^2)  \cap (\Omega_T \cup \Gamma)} u \le C \rho^\alpha
\end{equation}
with 
\begin{equation}
 \alpha = \min \left(- \log_{\sigma^2} ( 1 - 2^{-s}), \frac{n \kappa}{2} \right) \qquad C =  (\sfrac{\sigma^2}{\tilde{\rho}_0})^\alpha \max \left( \osc_{Q_{\tilde{\rho}_0}} u, 2^s \sigma^{\frac{n \kappa}{2}}  \tilde{\rho}_0^{\frac{n \kappa}{2}} \right)
\end{equation}
and any $s$ satisfying
\begin{equation}
s \ge  \left \lceil log_2 \frac{2 M }{\delta}\right \rceil\  +  (72 \omega_n \tilde{K}_P)^2  \gamma \theta \frac{2^{n+2 +2 \max \left(1, \frac{1+ \kappa}{r} \right) }}{\eta^2 (\sigma - 1)}
\end{equation}
\label{thm1}
\end{theo}
As a example of an application of the above mentioned theory, we present the proof of the result on the swirl of the axially symmetric Navier-Stokes flow in a cylinder. Before stating the result, let us introduce some quantities.\\
For $v_{r}, \; v_\phi,  \; v_z$ being the cylindrical components of three-dimensional vector field $u$ introduce quantity $u = r v_\phi$ called swirl. Let $v$ be a (weak) solution to Navier-Stokes system in a cylinder $\Omega_T$ with radius $R$:
\begin{equation}
\begin{aligned}
v,_t + v \cdot \nabla v - \nu \Delta v  = 0  \quad & \text{ in } \; \Omega_T\\
\text{div } u = 0  \quad & \text{ on } \; \Omega_T\\
v \cdot n = 0, \quad n \cdot {\mathbb{D}} (v) \cdot  \tau_i = 0  \quad & \text{ on } \; S_1^T \\
v \cdot n = 0  \quad & \text{ on } \; S_2^T \\
v_{|t=0} = v_0  \quad & \text{ in } \; \Omega \\
\end{aligned}
\label{ns}
\end{equation}
where $S_1$ denotes the curved part of boundary of the cylinder and $S_2$ - its (two-component) flat part. Consequently $u$ solves a following equation:
\begin{equation}
\begin{aligned}
u,_t + v \cdot \nabla u - \nu \Delta u + \nu \frac{u,_r}{r} = 0  \quad & \text{ in } \; \Omega_T\\
u,_r = \frac{2}{R} u  \quad & \text{ on } \; S_1^T \\
u \cdot n = 0  \quad & \text{ on } \; S_2^T \\
u_{|t=0} = u_0  \quad & \text{ in } \; \Omega \\
\end{aligned}
\label{swirl}
\end{equation}
\begin{theo}
Assume that $u \le M$ satisfies \eqref{swirl} with respective \eqref{ns} solution  $v \in L^{r'}(0,T; L^{q'} )$, $\frac{3}{q'} + \frac{2}{r'} = 1 - \frac{3}{2} \kappa $. Then $u \in B_N (\Omega_T,  M,  \gamma, r, \delta, \kappa)$ with any $\gamma \in (0, \nu), \delta \in \er$.
\label{swi1}
\end{theo}
As a corollary let us formulate
\begin{theo}
Assume for axially symmetric \eqref{ns} solution $v$ that $v_r, v_z \in L^{10} (\Omega_T)$, $r v_0$ is bounded and  in vicinity of the axis of symmetry $u_0$ is H\"older continuous with H\"older exponent $\frac{3}{2} \kappa $, $ \kappa \in ( 0, \frac{1}{3}]$. Then $u \in C^\alpha (\Omega^T) $.
\label{swi2}
\end{theo}
For the entire section, fix $x_0, t_0$ and supercylinder $Q_{(\min{\rho_0, 1}), 1}$ containing all further cylinders, where $\rho_0$ comes from the anti-cusp condition. Denote the boundary cylinder $Q(\sigma \rho, \theta (\sigma \rho)^2) $ by $Q_{\sigma \rho}$ and 
 \begin{equation}
\overline{m} \equiv \max_{Q_{\sigma \rho} \cap (\Omega_T \cup \Gamma) } u, \quad \underline{m} \equiv \min_{Q_{\sigma \rho} \cap (\Omega_T \cup \Gamma) } u, \quad \omega \equiv \osc_{Q_{\sigma \rho} \cap (\Omega_T \cup \Gamma) } u \;(= \overline{m} - \underline{m} )
\label{defwD}
 \end{equation}
\begin{lem}[Trichotomy for $B_D$]
Take $u \in B_D$. For any fixed  $\eta >0$ and $\sigma \in (1, 2]$ exists $s= s( \eta, \sigma)$ for which a following trichotomy holds for every time contraction parameter $\theta \le 1$:
\begin{equation}
\begin{aligned}
\text{either}  \quad (T1) &\qquad \omega \le 2^s \rho^{\min (\beta, \frac{n \kappa}{2})} \\
\text{or} \quad (T2) &\qquad\mu ( \{ (x,t) \in Q_\rho  \cap \Omega_T  : u(x,t) > \overline{m} - 2^{-s+1} \omega \} ) \le \eta \rho^{n+2}\\
\text{or} \quad (T2') &\qquad \mu ( \{ (x,t) \in Q_\rho  \cap \Omega_T : u(x,t) <  \underline{m}  + 2^{-s+1} \omega \} ) \le \eta \rho^{n+2}
\end{aligned}
\end{equation}
\label{lem73D}
\end{lem}
This lemma asserts quantitatively a following observation: for a function in $B_D$ either we (T1) control oscillations on $Q_{\sigma \rho}$  or  (T2) on a considerable fraction (in terms of Lebesgue measure) of a slightly smaller cylinder $Q_{ \rho}$ u is bounded away from its maximum or minimum in the bigger cylinder.
Define
\[
\lceil x \rceil = \inf_{\en \cup 0} \{c \ge x \}
\]
\begin{proof}
Assume that (T1) fails. 
Therefore 
\begin{equation}  
\frac{\omega}{2} > 2^{s-1} \rho^{\min (\beta, \frac{n \kappa}{2})} \ge  2^{s-1}\rho^\beta \ge c_\Gamma (\sigma \rho)^\beta \ge \osc_{Q (\sigma \rho, (\sigma \rho)^2) } u_{| \Gamma}
\label{tri1.2}
\end{equation}
where the last-but-one inequality is given by definition of $s$ i.e. \eqref{tri12} and last one by definition of $B_D$, point (iii).
Inequality \eqref{tri1.2} implies that 
\begin{equation}
\text{either } \; \max_{Q (\sigma \rho, (\sigma \rho)^2) } u_{| \Gamma} <  \overline{m} - \frac{\omega}{4}  \quad \text{ or } \; \min_{Q (\sigma \rho, (\sigma \rho)^2) } u_{| \Gamma} >  \underline{m} + \frac{\omega}{4} . 
\label{tri1.5}
\end{equation}
Assume that the former holds. Define
\begin{equation}
\begin{aligned}
r_0 &\equiv \left \lceil log_2 \frac{2 M}{ \delta} \right \rceil\\
k_r &\equiv \overline{m} - 2^{-r} \omega \quad \text{ for } r \ge r_0
\label{tri1.7}
\end{aligned}
\end{equation}
where $M, \delta$ are parameters of $B_D$.
Observe that \eqref{tri1.7} and assumption that the first possibility in \eqref{tri1.5} holds imply for $r_0 \ge 2$
\begin{equation}
k_r \ge \max \left( \max_{\Gamma \cap Q_{\sigma \rho}} u, \;  \esssup_{Q_{\sigma \rho} \cap (\Omega_T \cup \Gamma) } u  - \delta \right),
\label{tri2}
\end{equation}
so levels $k_r $ are admissible to \eqref{defB}. We show that (T2) is valid. For clarity the following main part of the proof is divided into a few steps
\begin{itemize}
\item[(i)]  Define a function in $Q_{ \rho}$
\begin{equation}
h(x,t) \equiv
\begin{cases}
\begin{aligned}
k_{r+1} - k_r &\quad \{ (x,t) \in Q_\rho  \cap (\Omega_T \cup \Gamma) : u (x,t) > k_{r+1}    \} \\
u (x,t) - k_r & \quad \{ (x,t) \in Q_\rho  \cap (\Omega_T \cup \Gamma) : k_r < u (x,t) \le k_{r+1}    \} \\
0 & \quad otherwise
\end{aligned}
\end{cases}
\label{tri2.5}
\end{equation}
Both $u \in B_D$ and \eqref{tri2} giving $u_{\Gamma \cap Q_{ \rho}} (x,t) - k_r \le 0$ imply that $h(t, \cdot) \in W^{1,1} ( B_{ \rho})$. Hence one can use Lemma \ref{gp}, choosing $\eta \equiv 1$
\begin{equation}
\int_{B_{ \rho}} h (t) \le K_P \rho^n \frac{ \mu^\frac{1}{n} (B_{ \rho})}{\mu_n (\{ x \in B_\rho : h (x, t) = 0 \} )} \int_{B_\rho} | \nabla h (t)|
\label{tri3}
\end{equation}
By definition $h = 0$ outside $\Omega$. Using this and Tchebytschev inequality one has from \eqref{tri3}
\begin{equation}
(k_{r+1} - k_r) \mu_n ( A_{k_{r+1}, \rho} (t) ) \le K_P  \frac{ {\omega_n}^{\frac{1}{n}} \rho^{n+1}}{\mu (B_{ \rho} \cap \Omega^c  )} \int_{ A_{k_{r}, \rho} (t) \backslash A_{k_{r+1}, \rho} (t) } | \nabla h (t)|
\label{tri3.5}
\end{equation}
where definition \eqref{defA} is used\footnote{ na raze wywalamy $\Gamma$ z T2 }.
In view of anti-cusp condition \eqref{ac}, \eqref{tri3.5} yields
\begin{equation}
\omega 2^{-(r+1)} \mu_n ( A_{\overline{m}  - \omega 2^{-(r+1)}, \rho} (t) ) \le \rho \frac{K_P  {\omega_n}^{\frac{1-n}{n}} }{\theta_0}\int_{ A_{k_{r}, \rho} (t) \backslash A_{k_{r+1}, \rho} (t) } | \nabla u^{(k_r)} (t)|
\label{tri4}
\end{equation}
Integrate \eqref{tri4} over $[ t_0 - \theta  \rho^2, t_0]$
\begin{multline}
\omega 2^{-(r+1)} \mu ( \{ (x,t) \in Q_\rho  \cap \Omega_T  : u(x,t) > \overline{m}  - 2^{-{(r+1)}} \omega \} ) \le \\
\rho \frac{K_P  {\omega_n}^{\frac{1-n}{n}} }{\theta_0} \int_{ t_0 - \theta  \rho^2}^{t_0} \int_{ A_{k_{r}, \rho} (t) \backslash A_{k_{r+1}, \rho} (t) } | \nabla u^{(k_r)}|
\end{multline}
Squaring this one has
\begin{multline}
\omega^2 4^{-(r+1)}  \mu^2 ( \{ (x,t) \in Q_\rho  \cap \Omega_T  : u(x,t) > \overline{m}  - 2^{-{(r+1)}} \omega \} ) \le \\  \rho^2 \frac{(K_P  {\omega_n}^{\frac{1-n}{n}} )^2}{\theta_0^2} \left[ \int_{ t_0 -\theta  \rho^2}^{t_0} \mu_n ( A_{k_{r}, \rho} (t) \backslash A_{k_{r+1}, \rho} (t) ) \right] \left[ \int_{ Q_\rho  \cap \Omega_T }  | \nabla u^{(k_r)}|^2 \right]
\label{tri5}
\end{multline}
\item[(ii)] To estimate term $\int_{ Q_\rho  \cap \Omega_T }  | \nabla u^{(k_r)}|^2 $ in \eqref{tri5} we use  the definition of $B_D$. Observe that \eqref{tri2} concludes that $k_r$ with $w = + u$ is admissible to \eqref{defB} in $Q_{\sigma \rho} $. This with
\begin{equation}
\xi (x,t)=
\begin{cases}
\begin{aligned}
 1 \quad& \text{in} \quad \overline{Q_{ \rho}} \\
0 \quad&  \text{outside} \quad \overline{Q_{ \sigma \rho}} \; \text{ and for } \; t = t_0 - \theta (\sigma \rho)^2 \\
\text{affine} \quad& \text{otherwise}
\end{aligned}
\end{cases}
\end{equation}
\begin{equation}
|\nabla \xi | \le |\rho (\sigma - 1)|^{-1}, \quad | \xi,_{t} | \le |\theta \rho^2 (\sigma^2 - 1)|^{-1}
\end{equation}
produces
\begin{multline}
 \gamma^{-1} \int_{Q_{\rho} \cap \Omega_T}| \nabla u^{(k_r)} |^2\le \\
( | \rho (\sigma - 1)|^{-2} +  |\theta \rho^2 (\sigma^2 - 1)|^{-1}) \int_{Q_{\sigma \rho} \cap \Omega_T}  | u^{(k_r)}|^2+\left(  \int_{t_0 - \theta (\sigma \rho)^2}^{t_0} \mu_n^\frac{r}{q} ({A_{k_r, \sigma \rho} (t) }) dt \right)^{\frac{2 (1+ \kappa)}{r}} 
\label{tri7}
\end{multline}
It holds
\[
 \int_{Q_{\sigma \rho} \cap \Omega_T}  | u^{(k_r)}|^2 =  \int_{Q_{\sigma \rho} \cap \Omega_T}  | (u - \overline{m}  + 2^{-r}\omega )^+|^2 \le 4^{-r}\omega^2 \omega_n \theta (\sigma \rho)^{n+2}
\]
\[
\left(  \int_{t_0 - \theta (\sigma \rho)^2}^{t_0} \mu_n^\frac{r}{q} ({A_{k_r, \sigma \rho} (t) }) dt \right)^{\frac{2 (1+ \kappa)}{r}}  \le \omega_n \left(  \theta (\sigma \rho)^{2 + \frac{n r}{q}}  \right)^{\frac{2 (1+ \kappa)}{r}} = \omega_n  \theta^\frac{2 (1+ \kappa)}{r} (\sigma \rho)^{ n (1+ \kappa)}
\]
where definition of $u^{(k_r)}$, $A_{k_r, \sigma \rho} (t) \subset B_{\sigma \rho}$ and $\frac{1}{r} + \frac{n}{2q} = \frac{n}{4} $ (see $B_D$ definition) are used. In view of the above two inequalities \eqref{tri7} implies
\begin{multline}
\int_{Q_{\rho} \cap \Omega_T}| \nabla u^{(k_r)} |^2 \le \\
 \gamma \omega_n \left[ ( | \rho (\sigma - 1)|^{-2} +  |\theta \rho^2 (\sigma^2 - 1)|^{-1})4^{-r}\omega^2 \theta (\sigma \rho)^{n+2}+ \theta^\frac{2 (1+ \kappa)}{r} (\sigma \rho)^{ n (1+ \kappa)} \right] \\
\le \gamma \omega_n \left[ 4^{-r}\omega^2 \theta (1 + \theta) \rho^n  \frac{\sigma^{n+2}}{\sigma - 1} + \theta^\frac{2 (1+ \kappa)}{r} (\sigma \rho)^{ n (1+ \kappa)} \right] 
\le \rho^n K_{\ref{tri7.5}}  \left[ 4^{-r}\omega^2+ \rho^{ n  \kappa} \right] .
\label{tri7.5}
\end{multline}
because by assumption $\theta \le 1$, one can take
\[
K_{\ref{tri7.5}} = \omega_n \gamma \frac{2^{n+2+2 \max \left(1, \frac{1+ \kappa}{r} \right) }}{\sigma - 1}  
\]
As $\rho^{ n  \kappa} = \rho^{\min ( 2 \beta, n  \kappa ) + (n  \kappa -  2 \beta)^+} \le \rho^{ (n  \kappa -  2 \beta)^+}  4^{-s}\omega^2$ holds by assumption that (T1) fails, one has for $r \in [ r_0, s] $ and $\rho \le 1$ 
\begin{equation}
\int_{Q_{\rho} \cap \Omega_T}| \nabla u^{(k_r)} |^2 \le K_{\ref{tri7.5}}  \rho^n   4^{-r}\omega^2
\label{tri8}
\end{equation}
\item[(iii)] Use \eqref{tri8} in \eqref{tri5} to get for  $r \in [ r_0, s] $
\begin{multline}
\omega^2 4^{-(r+1)}  \mu^2 ( \{ (x,t) \in Q_\rho  \cap \Omega_T  : u(x,t) > \overline{m}  - 2^{-{(r+1)}} \omega \} ) \le \\   4^{-r}\omega^2
 K_{\ref{lem73D}} \rho^{n+2}\left[ \int_{ t_0 -\theta  \rho^2}^{t_0} \mu_n ( A_{k_{r}, \rho} (t) \backslash A_{k_{r+1}, \rho} (t) ) \right]
 \label{tri9}
 \end{multline}
 for
 \begin{equation}
  K_{\ref{lem73D}} = \frac{(K_P {\omega_n}^{\frac{1-n}{n}} )^2}{\theta_0^2} K_{\ref{tri7.5}} = \left(\frac{1}{n}+ \omega_n \right)^2  \omega_n^{\frac{2- n}{n}} \gamma \frac{2^{3n+2+2 \max \left(1, \frac{1+ \kappa}{r} \right) }}{(\sigma - 1) \theta_0^2}  
\end{equation}
Divide \eqref{tri9} by $\omega^2 4^{-(r+1)} $; use $\overline{m}  - 2^{-{(r+1)}} \omega \le  \overline{m}  - 2^{-{(s-1)}} \omega$ for  $r \in [ r_0, s-2] $ and definition of $A_{k, \rho} (t)$
\begin{multline}
 \mu^2 ( \{ (x,t) \in Q_\rho  \cap \Omega_T  : u(x,t) > \overline{m}  - 2^{-{(s-1)}} \omega \} ) \le \\   4
 K_{\ref{lem73D}} \rho^{n+2} \mu ( \{ (x,t) \in Q_\rho  \cap \Omega_T  : u(x,t) \in (k_r, k_{r+1}]  \} ) .
  \label{tri10}
 \end{multline}
 To enable further control of constant, sum \eqref{tri10} over $r \in [ r_0, s-2] $
 \begin{multline}
 \mu^2 ( \{ (x,t) \in Q_\rho  \cap \Omega_T  : u(x,t) > \overline{m}  - 2^{-{(s-1)}} \omega \} ) \le \\  \frac{ 4
 K_{\ref{lem73D}}}{s -1 - r_0} \rho^{n+2} \mu ( \{ (x,t) \in Q_\rho  \cap \Omega_T  : u(x,t) \in (k_{r_0}, k_{s-1}]  \} ) \le  \frac{ 4 \omega_n K_{\ref{lem73D}} \theta}{s  - 1- r_0}   \rho^{2(n+2)} 
  \end{multline}
 So the main part of the proof results in
 \begin{equation}
  \mu ( \{ (x,t) \in Q_\rho  \cap \Omega_T  : u(x,t) > \overline{m}  - 2^{-{(s-1)}} \omega \} ) \le \sqrt{  \frac{ 4  \theta \omega_n K_{\ref{lem73D}} }{s -1- r_0}  } \rho^{(n+2)} 
  \label{tri11}
 \end{equation}
\end{itemize}
The proof concludes with a proper choice of $s$ satisfying:
\begin{equation}
 \sqrt{  \frac{ 4 \theta\omega_n K_{\ref{lem73D}} }{s  -1- r_0}  } \le \eta, \qquad 2 c_\Gamma \sigma^\beta   \le 2^s
 \label{tri12}
\end{equation}
The first inequality gives (T2) from \eqref{tri11} while the second allows for \eqref{tri1.2}. Recall \eqref{tri1.5}; its second alternative
is considered analogously as the above case, with $w = -u$ instead of $+u$, and yields (T2').
\end{proof}
Performing computation based on conditions \eqref{tri12} one has
\begin{rem}
In Lemma \ref{lem73D} any
\begin{equation}
s \ge 1+ \max \left[ \left \lceil log_2 \frac{2 M }{\delta}\right \rceil\  +  \theta \left(4 \left(\frac{1}{n}+ \omega_n \right)^2  \omega_n^{\frac{2}{n}} \gamma \frac{2^{3n+2+2 \max \left(1, \frac{1+ \kappa}{r} \right) }}{ \eta^2 (\sigma - 1) \theta_0^2} \right),  \; log_2 (2 c_\Gamma \sigma^\beta) \right]
\end{equation}
is admissible. One can choose $\theta$ small enough to shrink
\[
\theta \left(4 \left(\frac{1}{n}+ \omega_n \right)^2  \omega_n^{\frac{2}{n}} \gamma \frac{2^{3n+2+2 \max \left(1, \frac{1+ \kappa}{r} \right) }}{ \eta^2 (\sigma - 1) \theta_0^2} \right) 
\] 
as needed. Recalling that $\sigma \le 2$, sufficient condition for $s$ reads
\begin{equation}
s >   \left \lceil  log_2 \frac{2 M c_\Gamma }{\delta}   \right \rceil +2 + \beta
\end{equation}
\label{remsd}
\end{rem}
Below we state an analogous result to Lemma \ref{lem73N} for $B_N$. Recall that $ Q_{\sigma \rho} \equiv Q(\sigma \rho, \theta (\sigma \rho)^2) $. We define respective quantities without resorting to presently unknown boundary values
 \begin{equation}
\overline{m} \equiv \max_{Q_{\sigma \rho} \cap \Omega_T  } u, \quad \underline{m} \equiv \min_{Q_{\sigma \rho} \cap \Omega_T  } u, \quad \omega \equiv \osc_{Q_{\sigma \rho} \cap \Omega_T  } u \;(= \overline{m} - \underline{m} )
\label{defwN}
 \end{equation}
\begin{lem}[Trichotomy for $B_N$]
Assume that {\it anti-inner-cusp condition} \eqref{acn} holds. Take $u \in B_N$ and a cylinder $Q_{\rho \sigma}$ with the time contraction parameter satisfying
\begin{equation}
\theta \le \min \left(1, \frac{\theta_0}{2304 \gamma}, \; \left(  \frac{\omega^{-\frac{2(1+ \kappa) - q}{q}}_n}{128 \gamma} \right)^{\frac{r}{2 (1+ \kappa)}} \right) 
\label{trintheta}
\end{equation}
 Then for any fixed  $\eta >0$ and $\sigma \in (1, 2]$ exists $s= s(\eta,  \sigma)$, for which a following trichotomy holds:
 \begin{equation}
\begin{aligned}
\text{either}  \quad (T1) &\qquad \omega \le 2^s  \rho^{ \frac{n \kappa}{2}}  \\
\text{or} \quad (T2) &\qquad\mu ( \{ (x,t) \in Q_\rho  \cap \Omega_T  : u(x,t) > \overline{m} - 2^{-s+1} \omega \} ) \le \eta \rho^{n+2}\\
\text{or} \quad (T2') &\qquad \mu ( \{ (x,t) \in Q_\rho  \cap \Omega_T : u(x,t) <  \underline{m}  + 2^{-s+1} \omega \} ) \le \eta \rho^{n+2}
\end{aligned}
\end{equation}
\label{lem73N}
\end{lem}
An attempt to rewrite the proof of Lemma \eqref{lem73D} fails at obtaining \eqref{tri3.5} from \eqref{tri3}. Extrapolation of a truncated $u$ by zero outside $\Omega_T$, as in \eqref{tri2.5}, does not produce Sobolev function $h$ now, because boundary values of $u$ are not known. Thus one may extrapolate $u$ and define $h$ on $Q_\rho$ or restrict in definition of $h$ to $Q_\rho \cap \Omega_T$. In both cases we loose an easy way to control $ \mu_n(\{ h ( t) = 0 \}) $. Regaining this control poses  the main new point in the proof of Lemma \eqref{lem73N}. In the proof below we focus on this problem and sketch the part which overlaps with the previous proof.\\
Recall that
\[
\lceil x \rceil = \inf_{\en \cup 0} \{c \ge x \}
\]
\begin{proof}
Introduce 
\begin{equation}
\begin{aligned}
r_0 &\equiv \left \lceil log_2 \frac{2 M}{\delta}  \right \rceil, \\
k_r &\equiv \overline{m} - 2^{-r} \omega \quad \text{ for } r \ge r_0,
\label{trin5}
\end{aligned}
\end{equation}
levels $k_r$ are admissible to \eqref{defB}.
By definitions of $A, \omega, \overline{m}, \underline{m}  $ either
\begin{equation}
A^u_{\overline{m} - \frac{\omega}{2}, \rho} (t_0 - \theta \rho^2) \le \frac{1}{2} \mu_n (B_\rho \cap \Omega)
\label{trin1}
\end{equation}
or 
\begin{equation}
 A^{-u}_{\underline{m} - \frac{\omega}{2}, \rho} (t_0 - \theta \rho^2) \le \frac{1}{2} \mu_n (B_\rho \cap \Omega)
 \label{trin2}
\end{equation}
Consider case when \eqref{trin1} holds\footnote{The other one is performed analogously with $-u$ in place of $u$}, it implies for $r \ge 1$
\begin{equation}
A^u_{\overline{m} - \frac{\omega}{2^r}, \rho} (t_0 - \theta \rho^2) \le \frac{1}{2} \mu_n (B_\rho \cap \Omega).
\label{trin2.5}
\end{equation}
One can assume that both
\begin{equation}
 \max_{Q_\rho  \cap \Omega_T} u > \overline{m} - 2^{-s} \omega
 \label{trin3}
\end{equation}
holds, as otherwise (T2) holds with $\eta = 0$, and (T1) fails:
\begin{equation}
\omega > 2^s \rho^{ \frac{n \kappa}{2}} 
\label{trin4}
\end{equation}
The following essential part of the proof is divided into few steps.
\begin{itemize}
\item[(i)] 
Define 
\begin{equation}
h(x,t) \equiv
\begin{cases}
\begin{aligned}
k_{r+1} - k_r &\quad \{ (x,t) \in Q_\rho  \cap \Omega_T : u (x,t) > k_{r+1}    \} \\
u (x,t) - k_r & \quad \{ (x,t) \in Q_\rho  \cap \Omega_T  : k_r < u (x,t) \le k_{r+1}    \} \\
0 & \quad otherwise
\end{aligned}
\end{cases}
\label{trin6}
\end{equation}
As $h(t, \cdot) \in W^{1,1} ( B_{ \rho}  \cap \Omega)$, Lemma \ref{gpn} applies; choosing in it $\eta \equiv 1$ one has
\begin{equation}
\int_{B_{ \rho}  \cap \Omega} h (t) \le \tilde{K}_P  \frac{ \rho^{n+1}}{\mu_n (\{ x \in B_\rho \cap \Omega : h (x, t) = 0 \} )} \int_{B_\rho \cap \Omega} | \nabla h (t)|
\end{equation}
which implies
\begin{multline}
(k_{r+1} - k_r) \mu_n ( A^u_{k_{r+1}, \rho} (t) ) \le \\
  \tilde{K}_P \frac{  \rho^{n+1}}{\mu_n (\{ x \in B_\rho (x_0) \cap \Omega : u (x,t) \le k_r \})} \int_{ A^u_{k_{r}, \rho} (t) \backslash A^u_{k_{r+1}, \rho} (t) } | \nabla u^{(k_r)} (t)|
\label{trin7}
\end{multline}
As already remarked directly before the proof, we need  in \eqref{trin7} estimate of a following type
\begin{equation}
\mu_n (\{ x \in B_\rho (x_0) \cap \Omega : u (x,t) \le k_r \}) \ge \chi \omega_n \rho^n ;
\label{trin7.5}
\end{equation}
for some nonzero $\chi$. Such majorisation is done in the next step.
\item[(ii)] In \eqref{defB} take function $\eta (x)$ cutting off between $B_{\sfrac{\rho}{\lambda}} $ and $B_\rho$ with $\lambda > 1$
\begin{equation}
\eta (x,t)=
\begin{cases}
\begin{aligned}
 1 \quad& \text{in} \quad \overline{B_{\sfrac{\rho}{\lambda}}} \\
0 \quad&  \text{outside} \quad \overline{B_{  \rho}} \\
\text{affine} \quad& \text{otherwise}
\end{aligned}
\end{cases}
\end{equation}
\begin{equation}
|\nabla \eta | \le \frac{\lambda^2}{\rho^2 (\lambda - 1)^2}
\end{equation}
to obtain for $k_r$, admissible in view of \eqref{trin5},
\begin{multline}
\max_{t_0 - \theta \rho^2 \le t \le t_0} \int_{B_{\sfrac{\rho}{\lambda}} \cap \Omega} |u^{(k_r)} (t)|^2 \le \\\int_{B_{{\rho}} \cap \Omega} |u^{(k_r)} (t_0 - \theta \rho^2)|^2 + \gamma \left[ \frac{\lambda^2}{\rho^2 (\lambda - 1)^2} \int_{Q_{{\rho, \tau}} \cap \Omega} |u^{(k_r)} |^2 +  (\theta \rho^2)^{\frac{2 (1+ \kappa)}{r}}   \mu^{\frac{2 (1+ \kappa)}{q}}_n (B_{{\rho}} \cap \Omega) \right].
\label{trin15}
\end{multline}
Estimate the first summand of right hand-side of \eqref{trin15} by \eqref{trin2.5}; because left hand-side satisfies for $l > 0$
\[
\int_{B_{\sfrac{\rho}{\lambda}} \cap \Omega} |u^{(k_r)} (t)|^2 = \int_{A_{k_r, \sfrac{\rho}{\lambda} }(t)} |u^{(k_r)} (t)|^2 \ge \int_{A_{k_r + l, \sfrac{\rho}{\lambda}} (t)} |u^{(k_r)} (t)|^2 \ge l^2 \mu_n (A_{k_r + l, \sfrac{\rho}{\lambda}}(t))
\]
we have 
\begin{multline}
l^2 \mu_n (A_{k_r + l, \sfrac{\rho}{\lambda}}(t)) \le \frac{1}{2}  |u^{(k_r)} (t_0 - \theta \rho^2)|_{L^\infty(B_{{\rho}} \cap \Omega)}^2  \mu_n (B_\rho \cap \Omega) + \\
 \gamma \frac{\lambda^2}{\rho^2 (\lambda - 1)^2} \max_{t \in [t_0 - \theta \rho^2, t_0]} |u^{(k_r)}(t)|_{L^\infty(B_{{\rho}} \cap \Omega)}^2  \theta \rho^2 \mu_n (B_\rho \cap \Omega) +  \gamma (\theta \rho^2)^{\frac{2 (1+ \kappa)}{r}}   \mu^{\frac{2 (1+ \kappa)}{q}}_n (B_{{\rho}} \cap \Omega) 
\label{trin16}
\end{multline}
Define $H \equiv \max_{t \in [t_0 - \theta \rho^2, t_0]} |u^{(k_r)}(t)|_{L^\infty(B_{{\rho}} \cap \Omega)}$ and take $ l = H \psi$. Hence division  \eqref{trin16} by $l^2$ and estimate $\mu_n (B_\rho \cap \Omega) \le \mu_n (B_\rho)$ in its last summand yield
\begin{multline}
 \mu_n (A_{k_r + H \psi, \sfrac{\rho}{\lambda}}(t)) \\
  \le \frac{\mu_n (B_\rho \cap \Omega)}{\psi^2} \left[ \frac{1}{2}  
+ \gamma \left[ \frac{\lambda^2 \theta \rho^2}{\rho^2 (\lambda - 1)^2}  + \frac{ (\theta \rho^2)^{\frac{2 (1+ \kappa)}{r}}  \omega^\frac{2(1+ \kappa) - q}{q}_n \rho^{n\frac{2(1+\kappa) - q}{q}}}{H^2}  \right] \right]
\label{trin17}
\end{multline}
Assume 
\begin{equation}
r \le s-1
\label{trin18}
\end{equation}
this in tandem with definition of $H, k_r, \overline{m}$; \eqref{trin3}, $s \ge r_0$ gives
\begin{equation}
\begin{aligned}
H  \psi &= \left(\max_{Q_\rho  \cap \Omega_T} u -  [\overline{m} - 2^{-r} \omega ] \right) \psi \le  2^{-r} \omega \psi \\
H & = \max_{Q_\rho  \cap \Omega_T} u -  [\overline{m} - 2^{-r} \omega ] > 2^{-r} \omega - 2^{-s} \omega \ge  2^{-s} \omega \ge \rho^{ \frac{n \kappa}{2}} 
\end{aligned}
\label{trin19}
\end{equation}
The last inequality stems from  \eqref{trin4}. So \eqref{trin19}  in \eqref{trin17} yields
\begin{equation}
\begin{aligned}
\mu_n (A_{ \overline{m} - 2^{-r} \omega (1 -  \psi), \sfrac{\rho}{\lambda}}(t))  & \le \mu_n (A_{k_r + H \psi, \sfrac{\rho}{\lambda}}(t)) \\
\mu_n (A_{k_r + H \psi, \sfrac{\rho}{\lambda}}(t)) &\le \frac{\mu_n (B_\rho \cap \Omega)}{\psi^2} \left[ \frac{1}{2}  
+ \gamma \left[ \frac{\lambda^2 \theta}{ (\lambda - 1)^2}  +  \theta^{\frac{2 (1+ \kappa)}{r}}  \omega^\frac{2(1+ \kappa) - q}{q}_n \rho^{\nu}  \right] \right]
\end{aligned}
\end{equation}
where
\[
\nu = n\frac{2(1+\kappa) - q}{q} + {\frac{4 (1+ \kappa)}{r}}  - n \kappa = 4 ( 1+  \kappa) \left[ \frac{1}{r} + \frac{n}{2q} \right] - n ( 1+  \kappa) = 0
\]
because by $B_N$ definition $ \frac{1}{r} + \frac{n}{2q} = \frac{n}{4}$, hence for $\psi = \sfrac{3}{4}$
\begin{equation}
\mu_n (A_{\overline{m} - 2^{-(r+2)} \omega, \sfrac{\rho}{\lambda}}(t)) \le \mu_n (B_\rho \cap \Omega) \left[ \frac{8}{9}  
+ \frac{16 \gamma }{9} \left[ \frac{\lambda^2 \theta}{ (\lambda - 1)^2}  +  \theta^{\frac{2 (1+ \kappa)}{r}}  \omega^\frac{2(1+ \kappa) - q}{q}_n  \right] \right]
\label{trin20}
\end{equation}
Combine \eqref{trin20} with
\begin{equation}
\mu_n (A_{k, \rho}(t)) \le \mu_n (A_{k, \sfrac{\rho}{\lambda}}(t)) + \omega_n \left( \frac{\lambda - 1}{\lambda} \right)^n \rho^n \le \mu_n (A_{k, \sfrac{\rho}{\lambda}}(t)) + \frac{1}{18}  \mu_n (B_\rho \cap \Omega),
\end{equation}
where the equality holds for 
\begin{equation}
\left( \frac{\lambda - 1}{\lambda} \right)^n =  \frac{\theta_0}{18} 
\label{trin20.5}
\end{equation}
thanks to {\it anti-inner-cusp condition} \eqref{acn}, to get
\begin{equation}
\mu_n (A_{\overline{m} - 2^{-(r+2)} \omega, \rho}(t)) \le \mu_n (B_\rho \cap \Omega)\left[ \frac{17}{18}  
+ \frac{16 \gamma}{9} \left[ \frac{\lambda^2 \theta}{ (\lambda - 1)^2}  +  \theta^{\frac{2 (1+ \kappa)}{r}}  \omega^\frac{2(1+ \kappa) - q}{q}_n  \right] \right]
\label{trin21}
\end{equation}
Recall that estimating \eqref{trin7.5} is the aim of this step of the proof; \eqref{trin21} implies\begin{multline}
\mu_n (\{ x \in B_\rho (x_0) \cap \Omega : u (x,t) \le k_{r+2}  \}) \ge \\
 \mu_n (B_\rho \cap \Omega)\left[ \frac{1}{18}  
- \frac{16 \gamma}{9} \left[ \frac{\lambda^2 \theta}{ (\lambda - 1)^2}  +  \theta^{\frac{2 (1+ \kappa)}{r}}  \omega^\frac{2(1+ \kappa) - q}{q}_n  \right] \right] \ge \\
\theta_0 \omega_n \rho^n \left[ \frac{1}{18}  
- \frac{16 \gamma}{9} \left[  \theta \left(\frac{18}{\theta_0} \right)^{\frac{2}{n}}  +  \theta^{\frac{2 (1+ \kappa)}{r}}  \omega^\frac{2(1+ \kappa) - q}{q}_n  \right] \right],
\label{trin22}
\end{multline}
where the last inequality results from  {\it anti-inner-cusp condition} \eqref{acn} and  \eqref{trin20.5}. Therefore for validity of  \eqref{trin7.5} one needs
\begin{equation}
 \chi \equiv 
\theta_0 \left[ \frac{1}{18}  
- \frac{16 \gamma}{9} \left[  \theta \left(\frac{18}{\theta_0} \right)^{\frac{2}{n}}  +  \theta^{\frac{2 (1+ \kappa)}{r}}  \omega^\frac{2(1+ \kappa) - q}{q}_n  \right] \right] > 0.
\label{trin23}
\end{equation}
For any $\theta$ satisfying \eqref{trintheta} one computes $\chi \ge \sfrac{1}{36}$. Summing up, estimate \eqref{trin7.5} indeed holds in a following form
\begin{equation}
\mu_n (\{ x \in B_\rho (x_0) \cap \Omega : u (x,t) \le k_{r+2} \}) \ge \frac{ \omega_n \rho^n}{36}
\label{trin24}
\end{equation}
provided \eqref{trin18} holds, i.e $r \le s-1$.
\item[(iii)]  Use of \eqref{trin24} in \eqref{trin7} gives for $ {t \in [t_0 - \theta \rho^2, t_0]} $
\begin{equation}
(k_{r+1} - k_r) \mu_n ( A^u_{k_{r+1}, \rho} (t) ) \le \frac{36} {\omega_n}  \tilde{K}_P \rho \int_{ A^u_{k_{r}, \rho} (t) \backslash A^u_{k_{r+1}, \rho} (t) } | \nabla u^{(k_r)} (t)|
\label{trin25}
\end{equation}
which is an exact analogue of \eqref{tri4} in the proof of Lemma \ref{lem73D}. Let us only sketch the remainder of the proof, as from \eqref{trin25} it progresses along the lines of Lemma \ref{lem73D}. Estimate \eqref{trin25} implies an analogue of \eqref{tri5}
\begin{multline}
\omega^2 4^{-(r+1)}  \mu^2 ( \{ (x,t) \in Q_\rho  \cap \Omega_T  : u(x,t) > \overline{m}  - 2^{-{r+1}} \omega \} ) \le \\  \rho^2 \left( \frac{36} {\omega_n}  \tilde{K}_P \right)^2 \left[ \int_{ t_0 -\theta  \rho^2}^{t_0} \mu_n ( A_{k_{r}, \rho} (t) \backslash A_{k_{r+1}, \rho} (t) ) \right] \left[ \int_{ Q_\rho  \cap \Omega_T }  | \nabla u^{(k_r)}|^2 \right]
 \label{trin26}
\end{multline}
This, combined with computation rewritten from point (ii) of Lemma \ref{lem73D}, implies for $r \in [ r_0, s] $, a  \eqref{tri9} counterpart
\begin{multline}
\omega^2 4^{-(r+1)}  \mu^2 ( \{ (x,t) \in Q_\rho  \cap \Omega_T  : u(x,t) > \overline{m}  - 2^{-{r+1}} \omega \} ) \le \\   4^{-r}\omega^2
 K_{\ref{lem73N}} \rho^{n+2}\left[ \int_{ t_0 -\theta  \rho^2}^{t_0} \mu_n ( A_{k_{r}, \rho} (t) \backslash A_{k_{r+1}, \rho} (t) ) \right]
 \label{trin27}
 \end{multline}
 with
 \begin{equation}
 K_{\ref{lem73N}} =(36  \tilde{K}_P)^2  \omega_n \gamma \frac{2^{n+2 +2 \max \left(1, \frac{1+ \kappa}{r} \right) }}{\sigma - 1}  
 \end{equation}
As in the $B_D$ case, summing \eqref{trin27} over $r \in [ r_0, s-2] $ gives
 \begin{equation}
  \mu ( \{ (x,t) \in Q_\rho  \cap \Omega_T  : u(x,t) > \overline{m}  - 2^{-{(s-1)}} \omega \} ) \le \sqrt{  \frac{ 4 \omega_n \theta K_{\ref{lem73N}} }{s - 1- r_0}  } \rho^{(n+2)}
  \label{trin28}
 \end{equation}
\end{itemize}
Estimate \eqref{trin28} implies (T2) provided
 \begin{equation}
 \sqrt{  \frac{ 4 \theta \omega_n K_{\ref{lem73N}} }{s - 1- r_0}  } \rho^{(n+2)} \le \eta
  \label{trin29}
 \end{equation}
\end{proof}
As before, formulate
\begin{rem}
In Lemma \ref{lem73N} any
\begin{equation}
s \ge  \left \lceil log_2 \frac{2 M }{\delta}\right \rceil\  +  (72 \omega_n \tilde{K}_P)^2  \gamma \theta \frac{2^{n+2 +2 \max \left(1, \frac{1+ \kappa}{r} \right) }}{\eta^2 (\sigma - 1)}
\end{equation}
is good. Therefore, taking $\theta$ such that
\[
 (72 \omega_n \tilde{K}_P)^2  \gamma \theta \frac{2^{n+2 +2 \max \left(1, \frac{1+ \kappa}{r} \right) }}{\eta^2 (\sigma - 1)} \le 1
 \]
 we see that any
 \begin{equation}
 s \ge \max \left( 3,  \left \lceil log_2 \frac{2 M }{\delta}\right \rceil\  + 1 \right)
 \end{equation}
 is admissible.
 \label{remsn}
\end{rem}
In the following two lemmas we elaborate the case, when alternative (T2), (T2') of Lemmas \ref{lem73D}, \ref{lem73N} hold.
\begin{lem}[Vanishing measure]
Take $u \in B_D$ or $B_N$. There exists $\eta_0 > 0$ such that for any level $k$ admissible to \eqref{defB} and boundary cylinder $Q_\rho \subset \Omega_T$ inequality
\begin{equation}
\mu ( \{ (x,t) \in Q_\rho  \cap \Omega_T  : u(x,t) > k  \} ) \le \eta_0 \rho^{n+2} 
\end{equation}
implies that
\begin{equation}
\text{either } \max_{Q_\rho} u - k  < \rho^{ \frac{n \kappa}{2}} \quad \text{ or } \quad \mu \left( \left \{ (x,t) \in Q_{\sfrac{\rho}{\sigma}}  \cap \Omega_T  : u(x,t) > \frac{k  + \max_{Q_\rho} u}{2} \right \} \right) = 0
\end{equation}
\label{72D}
\end{lem}
\begin{lem}
\label{72N}
\end{lem}
Currently we are ready to present key results, enabling quantitative control of oscillations.
\begin{lem}[Oscillation control for $B_D$]
Assume that {\it anti-outer-cusp condition} \eqref{ac} holds. Take a cylinder $Q_{\rho \sigma}$ with the time contraction parameter $\theta \le 1$ and $u \in B_D$.  Then for any fixed $\sigma \in (1, 2]$ exists $s= s(\sigma)$, for which either
\begin{equation}
\begin{aligned}
&\mathrm{(O1) } \qquad \osc_{Q_{\sfrac{\rho}{\sigma}} \cap \Omega_T  } u \le 2^s \rho^{ \frac{n \kappa}{2}} \\
\text{or}&\\
&\mathrm{(O2)} \qquad \osc_{Q_{\sfrac{\rho}{\sigma}} \cap \Omega_T  } u \le (1 - 2^{-s}) \osc_{Q_{\rho\sigma} \cap \Omega_T  } u
\end{aligned}
\end{equation}
for every $\rho \le \frac{\rho_0}{\sigma}$, where $\rho_0$ stems from {\it anti-outer-cusp condition} \eqref{ac}.
To be precise, $s = s(\eta_0, \sigma)$ from Lemma \ref{lem73D} with $\eta_0$ fixed by Lemma \ref{72D}.
\label{oscD}
\end{lem}
\begin{lem}[Oscillation control for $B_N$]

Assume that {\it anti-inner-cusp condition} \eqref{acn} holds. Take $u \in B_N$ and a cylinder $Q_{\rho \sigma}$ with the time contraction parameter satisfying \eqref{trintheta}.
 Then for any fixed $\sigma \in (1, 2]$ exists $s= s(\sigma)$, for which either
\begin{equation}
\begin{aligned}
&\mathrm{(O1) } \qquad \osc_{Q_{\sfrac{\rho}{\sigma}} \cap \Omega_T  } u \le 2^s \rho^{ \frac{n \kappa}{2}} \\
\text{or}&\\
&\mathrm{(O2)} \qquad \osc_{Q_{\sfrac{\rho}{\sigma}} \cap \Omega_T  } u \le (1 - 2^{-s}) \osc_{Q_{\rho\sigma} \cap \Omega_T  } u
\end{aligned}
\end{equation}
for every $\rho \le \frac{\rho_0}{\sigma}$, where $\rho_0$ stems from {\it anti-outer-cusp condition} \eqref{acn}.
To be precise, $s = s(\eta_0, \sigma)$ from Lemma \ref{lem73N} with $\eta_0$ fixed by Lemma \ref{72N}.
\label{oscN}
\end{lem}
Having Lemmas \ref{lem73N}, \ref{72N} and \ref{lem73D}, \ref{72D} one shows Lemmas \ref{oscN}, \ref{oscD} as Lemma 7.4 from \cite{lsu}, Chapter II.7 (for some more details in case of cylinders scaled by factor $\sigma$ instead of $2$ as in \cite{lsu}, compare Lemma 5.3 from \cite{zaj}). Because the argument is straightforward, we present it below for reader's convenience.
\begin{proof}[Proof of Lemmas \ref{oscD}, \ref{oscN}] As we consider both $B_N$ and $B_D$ case, for convenience we refer to either Lemma \ref{lem73D} or \ref{lem73N} as {\it trichotomy lemma} and to Lemma \ref{72D} or \ref{72N} as {\it vanishing measure lemma}. Fix $\sigma_0 \in (1,2]$. Set the smallness parameter $\eta$ equal to $\eta_0$ from {\it vanishing measure lemmas}. Take as  $s_0 = s (\eta_0, \sigma_0)$ from {\it trichotomy lemmas}. Suppose (O1) does not hold. Therefore, in view of {\it trichotomy lemmas}, (T2) or (T2') is valid. Focus on the case when (T2) holds\footnote{ again, the other one is performed in the same way, with $-u$ instead of $u$}:
\[
\mu ( \{ (x,t) \in Q_\rho  \cap \Omega_T  : u(x,t) > k_{s_0 -1}  \} ) \le \eta_0 \rho^{n+2} \text{ \quad with \quad }
k_{s_0 -1} = \overline{m} - 2^{-(s_0-1)} \omega
\]
In view of definition of $s_0$, level $k_{s_0 -1}$ is admissible in \eqref{defB}. Observe that for a fixed $\eta_0$ one has thesis of {\it vanishing measure lemmas} for every level $k$ admissible to \eqref{defB} independently from $s$, so there is no loop in the above choice of parameters.
This allows us via {\it vanishing measure lemmas} to state 
that either
\begin{equation}
\max_{Q_\rho} u - k_{s_0 -1}  < \rho^{ \frac{n \kappa}{2}}
\label{osc5}
\end{equation}
or
\begin{equation}
\mu \left( \left \{ (x,t) \in Q_{\sfrac{\rho}{\sigma}}  \cap \Omega_T  : u(x,t) > \frac{k_{s_0 -1}  + \max_{Q_\rho} u}{2} \right \} \right) = 0
\label{osc6}
\end{equation}
In view of definition of $ k_{s_0 -1} $ and assumption that (O1) fails, \eqref{osc5} yields
\begin{equation}
\max_{Q_{\sfrac{\rho}{\sigma}} \cap \Omega_T  } u  \le \max_{Q_\rho \cap \Omega_T } u < k_{s_0 -1}  + \rho^{ \frac{n \kappa}{2}} <  \overline{m} - 2^{-(s_0-1)} \omega + 2^{-s_0}  \osc_{Q_{\sfrac{\rho}{\sigma}} \cap \Omega_T  } u \le \overline{m}  - 2^{-s_0} \omega.
\label{osc7}
\end{equation}
In case of validity of \eqref{osc6} holds
\begin{equation}
 \max_{Q_{\sfrac{\rho}{\sigma}} \cap \Omega_T  } u \le \frac{ \overline{m} - 2^{-(s_0-1)} \omega +\max_{Q_\rho} u}{2} \le \overline{m} - 2^{-s_0} \omega
 \label{osc8}
\end{equation}
Therefore \eqref{osc7}, \eqref{osc8} imply thesis, because
\begin{equation}
\osc_{Q_{\sfrac{\rho}{\sigma}} \cap \Omega_T  } u =  \max_{Q_{\sfrac{\rho}{\sigma}} \cap \Omega_T  } u - \min_{Q_{\sfrac{\rho}{\sigma}} \cap \Omega_T  } u < \overline{m}  - 2^{-s_0} \omega - \min_{Q_{\sfrac{\rho}{\sigma}} \cap \Omega_T  } u \le (1 - 2^{-s_0}) \osc_{Q_{\rho\sigma} \cap \Omega_T  } u.
\end{equation}
\end{proof}
We are ready to derive H\"older regularity result from the above formulated {\it oscillation control lemmas}, i.e. Lemmas \ref{oscD}, \ref{oscN}, use the fact below, where as usual $Q_{\rho}$ denotes $ Q(\rho, \theta \rho^2) $
\begin{lem}
Fix $\theta$.   If measurable,  bounded $u :  Q_{\rho_0} \cap \Omega_T \rightarrow \er$ satisfies for $\eta < 1$, $b >1$
\begin{equation}
 \text{either } \osc_{Q_{{\rho}} \cap \Omega_T  } u \le \eta\osc_{Q_{b \rho} \cap \Omega_T  } u \quad \text{or} \osc_{Q_{{\rho}} \cap \Omega_T  } u \le c_1 \rho^\delta
 \end{equation}
then $u$ is H\"older continuous; more precisely for every $\rho \le b^{-1} \rho_0$ holds
\[
 \osc_{Q_{\rho} \cap \Omega_T } u \le C \rho^\alpha
 \]
with 
 \[
 \alpha = \min (- \log_{b} \eta, \delta) \quad C = (\sfrac{b}{\rho_0})^\alpha \max \left( \osc_{Q_{\rho_0}} u, c_1 \rho_0^\delta \right)
 \]
 \label{hold}
 \end{lem}
 The proof of the above Lemma \ref{hold} con be found in \cite{lsu}, Chapter II.5. Finally we formulate
 \begin{proof}[Proof of Theorems \ref{thm}, \ref{thm1}]
 Combination of {\it oscillation control lemmas}, i.e. Lemmas  \ref{oscD}, \ref{oscN} and Lemma \ref{hold} gives the main result. Exact estimates for quantity $s$ is given by Remarks  \ref{remsd},  \ref{remsn}.
 \end{proof}
  \begin{proof}[Proof of Theorems \ref{swi1}, \ref{swi2}]
  Having theory on boundary DeGiorgi-Ladyzhenskaya classes, we perform these proofs exactly as proofs of Lemma 3.2 and Main Theorem in \cite{zaj2}.
 \end{proof}



{\bf Bibliography}

\begin{thebibliography}{99}
 

\bibitem{lsu} Ladyzhenskaya, O. A., Solonnikov V. A., Uraltseva, N. N. {\it Linear and quasilinear equations of parabolic type},
\bibitem{zaj} Zajaczkowski, W. M. {\it Global regular axially symmetric solutions to the Navier-Stokes equations in a periodic cylinder},
\bibitem{zaj2} Zajaczkowski, W. M. {\it The H\"older continuity of the swirl for the Navier-Stokes motions},



 \end{thebibliography}
\end{document}